\documentclass[a4paper,fleqn]{cas-dc}

\usepackage[utf8]{inputenc}
\usepackage[T1]{fontenc}   

\usepackage[authoryear]{natbib}

\usepackage{dsfont} 
\usepackage{subcaption}
\usepackage{stmaryrd}

\def\tsc#1{\csdef{#1}{\textsc{\lowercase{#1}}\xspace}}
\tsc{WGM}
\tsc{QE}
\newcommand{\esp}[2]{\mathds{E}_{#1}\left[#2\right]}

\begin{document}
\let\WriteBookmarks\relax
\def\floatpagepagefraction{1}
\def\textpagefraction{.001}

\shorttitle{A metamodel for confined yield stress flows and parameters' estimation}

\shortauthors{C.Berger et al.}  

\title[mode = title]{A metamodel for confined yield stress flows and parameters' estimation}  

\tnotemark[1] 

\tnotetext[1]{Supplementary code on \href{https://zenodo.org/records/8377205}{Zenodo}, ref
  8377205. This is the authors' version as of: February 23, 2024. Now accepted in:
  Rheologica Acta (2024)
  \href{https://doi.org/10.1007/s00397-024-01436-0}{https://doi.org/10.1007/s00397-024-01436-0}}

%

\author[1]{Clément Berger}[type=editor,
auid=000,
bioid=1,
orcid=0009-0003-0055-0803,
]


\fnmark[1]

\fntext[1]{Funded by a CDSN grant from the French MESRI} 

\ead{clement.berger@ens-lyon.fr}

\credit{Conceptualization, Methodology, Software, Writing}

\affiliation[1]{organization={UMPA CNRS UMR 5669 ENS de Lyon},
  addressline={46 Allée d'Italie}, 
  city={Lyon},
  postcode={69364}, 
  state={Rhône},
  country={France}}

\author[1]{David Coulette}[auid=001,
bioid=2,
orcid=0000-0002-3556-0089,
]

\fnmark[2]

\ead{david.coulette@ens-lyon.fr}

\credit{Data Curation, Software}

\author[1,2]{Paul Vigneaux}[type=editor,
auid=003, 
bioid=4, 
orcid=0000-0001-6606-5446,
]

\fnmark[2]

\cormark[1]

\ead{paul.vigneaux@math.cnrs.fr}


\credit{Conceptualization, Methodology, Software, Writing, Funding}

\affiliation[2]{organization={LAMFA CNRS UMR 7352 Université de Picardie Jules Verne},
  addressline={33 Rue Saint Leu}, 
  city={Amiens},
  postcode={80039}, 
  state={Somme},
  country={France}}

\cortext[1]{Corresponding author}

\fntext[2]{Funded by the French ANR under grant number ANR-20-CE46-0006}

\begin{abstract}
  With the growing demand of mineral consumption, the management of the mining waste is
  crucial. Cemented paste backfill (CPB) is one of the techniques developed by the mining
  industry to fill the voids generated by the excavation of underground spaces. The
  CPB process is the subject of various studies aimed at optimizing its implementation in the
  field. In this article, we focus on the modelling of the backfill phase where it has
  been shown in [Vigneaux et al., Cem. Concr. Res. 164 (2023) 107038] that a viscoplastic
  lubrication model can be used to describe CPB experiments. The aim here is to propose
  an accelerated method for performing the parameters' estimation of the properties of the
  paste (typically its rheological properties), with an inverse problem procedure based on
  observed height profiles of the paste. The inversion procedure is based on a
  metamodel built from an initial partial differential equation model, thanks to a
  Polynomial Chaos Expansion coupled with a Principal Component Analysis. 
\end{abstract}


\maketitle

\section{Introduction}\label{intro}

In this article, we focus on a fast model for the simulation and parameters' estimation of
thin sheets of viscoplastic material aiming at filling elongated cavities. A typical
application would be the cemented paste backfill in underground stopes.

With the growing demand of mineral consumption, the management of the mining waste is
crucial. Cemented paste backfill (CPB) is one of the techniques developed by the mining
industry to fill the voids generated by the excavation of underground spaces. CPB has been
documented in Germany in the 1980's \cite{landriault_they_2006}. It is now used in a growing
number of countries across the world, including Australia, Canada and China
\cite{Yilmaz2017paste}. The paste is obtained by mixing tailings coming from the
extraction, with water and a hydraulic binder. Mixtures are then transported from the
surface plant to the underground openings through pipes. Addition of the binder (cement
like) is crucial for the final strength and stability of the backfill. Among the
advantages of the CPB, one expects:
\begin{itemize}
\item the ability of reusing the mining waste by reintroducing it in the mine. This would
  prevent to use surface spaces on the ground (e.g., tailings dams), with all the potential
  environmental impacts this can have (\cite{UNO92} and \cite{AzLi10}).
\item By filling the voids in underground spaces with these pastes, the ground is thus
  consolidated, leading to two advantages: (i) the risk of subsidence
  \cite{keller_introduction_2008} is lowered and (ii) there is no need to leave unexploited
  pillars as done originally for stability, since the whole cavity can then be backfilled
  with the cemented paste (like in the cut-and-fill mining methodology). A better
  exploitation of the mining deposit can thus be expected.
\end{itemize} The practical implementation of CPB is a complex process where optimization
needs to be done on several aspects: the composition of the CPB (chemical and mineralogical
characteristics) and its resulting physical properties in relation with the transportation
phase (where fluidity would ease the pumping) and finally the consolidating phase (where
high mechanical strength is needed for an efficient stability of the resulting backfilled
volume) \cite{Yilmaz2017paste}.

In the present paper, we focus more specifically on the modelling and simulation of the
filling phase where we consider the slurry of a viscoplastic material in a \emph{bounded
  domain}. We refer to \cite{VSF21} for a more detailed description on the modelling where
it was shown that a viscoplastic lubrication model successfully described laboratory
experiments of such slurries made of typical CPB material. Using this model, we present
here a new methodology to perform parameters' estimation of this model thanks to an
inverse problem using \emph{observations} of the surface height of the slurry. The
objective of the present article is thus to provide a fast algorithm which allows to simulate
an observed paste flow and to determine its viscoplastic properties, typically the viscosity
(or consistency), the yield stress and the power index of the constitutive law of the
associated rheology (namely, the Herschel-Bulkley law). As a matter of fact the rheology
is one of the most important properties of CPB material \cite{Yilmaz2017paste,RC05}. Having
these parameters, we can then solve the direct problem to simulate further in time and
optimize the filling process of the cavity.

To perform the inverse problem in a faster way, we develop a so called \emph{metamodel}
which is a \emph{precise} approximation of the lubrication model (in its original partial
differential equation (PDE) form) but \emph{much faster} to evaluate than solving the
PDE. (Metamodels are also sometimes called \emph{surrogate} models.) Indeed, when
considering the closed stope configuration, there is no explicit solution for the height
profile and a numerical resolution is required. The design of this algorithm borrows from
previous works \cite{PCA_PCE_1,PCA_PCE_2} but seems to have been rarely adapted in the
context of PDEs. We develop here a specific and complete framework for the model of
\cite{VSF21}. Of note, we also provide in an open access repository the data and code of
the metamodel
\footnote{\href{https://zenodo.org/records/8377205}{https://zenodo.org/records/8377205}},
so that readers can experiment and identify the parameters of their own CPB lab experiments.

The paper is organized as follows. In section \ref{sec:dipr}, we detail the model under
consideration and the numerical resolution of the direct problem, as well as the links
with the inverse problem. We then present in section \ref{sec:metam} the construction of
the metamodel, which mixes Polynomial Chaos Expansion (PCE) and a Principal Component
Analysis (PCA). In section \ref{sec:resul}, we test the metamodel performances depending
on the characteristics of the PCA and the PCE. We also measure the ability to estimate the
parameters when the metamodel is coupled to a Nelder-Mead algorithm to solve the inverse
problem for both synthetic data (\emph{in silico}) or for the laboratory experiments data
of \cite{VSF21}.

\section{Formulation of the direct problem}\label{sec:dipr}

\subsection{Original PDE model}\label{sec:descripmod}

The starting model is a so-called lubrication model for yield stress flows in confined
geometry. It is borrowed from \cite{BCRS06} developed in unconfined geometry (see also
earlier works \cite{LiuMei89,CPA96,HG1998}). It was then studied in the case of the flow
in a confined cavity by \cite{VSF21}. We refer to \cite{BCRS06} and \cite{VSF21} for more
details.

\begin{figure}[h!]
  \centering
  \includegraphics[width=0.4\textwidth]{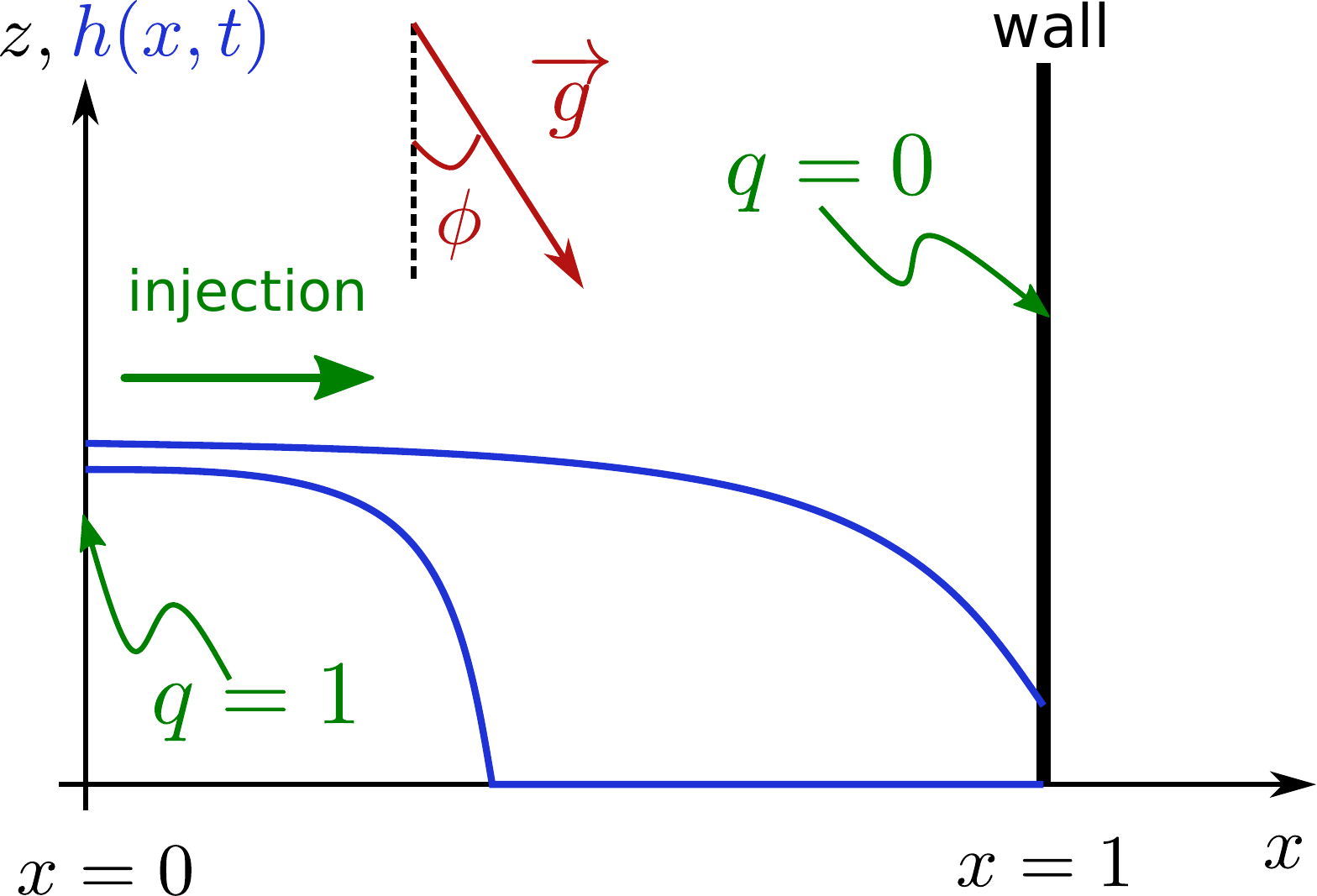}
  \caption{1D model, non-dimensionalized variables. The blue curve is the height $h$ of
    the material at two successive times (the flow is from left to right). Note that the
    gravity vector is inclined with an angle $\phi$ since the $x$ axis is assumed to be
    inclined downslope from the horizontal.}\label{fig:geom}
\end{figure}

We consider a closed cavity in 1D, inclined at an angle $\phi$ to the
horizontal (see figure \ref{fig:geom}). In non-dimensional variables the length of this
cavity is equal to $1$. We assume that the material follows a Herschel-Bulkley
constitutive law (whose characteristic parameters: yield stress, consistency and power
index will be described later). It is assumed that the injected flow rate of the material
is imposed on the left of the domain. The non-dimensionalization process leads to $q=1$ as
the left boundary condition for the model at $x=0$. While on the right, at $x=1$, a wall
condition is imposed with $q=0$. The unknown of the problem that is solved is the
dimensionless height $h$ of the material at time $t$ and distance $x$ from the injection
point. Note that we assume that there is no material at $t=0$ (the cavity is empty), so
that the initial condition of the PDE is : $h(x,t=0)=0$ for all $x \in [0,1]$. The
evolution of $h$ is the solution of the following partial differential equation (PDE):
\begin{equation}
  \label{eq:lub1D}
  \frac{\partial }{\partial t} h(x,t) + \frac{\partial }{\partial x}  q(x,t) = 0,
  \forall x\in ]0,1[, \forall t>0,
\end{equation}
where $q(x,t)$ is the flux function defined as:
\begin{align}
  \label{eq:qdef}
  q(x,t) & = \mbox{sgn} \left ( S - \frac{\partial h(x,t)}{\partial x} \right ) \;
           \left | S - \frac{\partial h(x,t)}{\partial x} \right |^{1/n}  \\
         & ~~~ \times  \frac{n {Y(x,t)}^{1+1/n} }{(n+1)(2n+1)} \left ( (2n+1)h(x,t) - nY(x,t) \right ), \nonumber
\end{align}
where $\mbox{sgn}(.)$ is the sign function and $n$ (dimensionless, typically in
$[0.3;1.2]$) is the power index of the aforementioned constitutive law. The variable
$Y(x,t)$ encodes the yield surface (i.e., the height above which the material behaves like
a pseudo-plug) and is given by:
\begin{equation}
  \label{eq:ypp}
  Y(x,t) = \max \left ( h(x,t) - \frac B {\left | S - \frac{\partial h}{\partial x}(x,t)
\right |}\, , \, 0 \right ),
\end{equation}
where $B$ is the \emph{Bingham} number and $S$ the so-called \emph{slope} parameter, in
this context.

Note that, as such, the independent parameters of this model are $B$, $S$ and $n$. Knowing
$(B,S,n)$, using aforementioned initial condition on $h$ and boundary conditions on $q$,
one can compute the evolution of $h$ in the cavity. We now detail these parameters in
terms of the associated dimensional variables. The slope parameter $S$ is given by:
\begin{equation}
  \label{eq:defs}
  S = \frac{\tan \phi}{\epsilon},
\end{equation}
where $\epsilon = \hat{H_0}/\hat{L}$ is the aspect ratio between the typical height of the
material ($\hat{H_0}$ in $m$) and the typical length of the cavity ($\hat{L}$ in $m$, in
the $x$ direction depicted in figure \ref{fig:geom}). In the lubrication theory (thin
film), it is assumed that $\epsilon \ll 1$.

The Bingham number is defined as:
\begin{equation}
  \label{eq:defbi}
  B = \frac{\hat{H}_0 \hat{\tau}_y}{\hat{\rho} \hat{\nu} \hat{U}},
\end{equation}
where $\hat{\tau}_y$ ($kg.m^{-1}.s^{-2}$ or $Pa$) is the yield stress of the viscoplastic
material, $\hat{\rho}$ is the density ($kg.m^{-3}$), $\hat{U}$ is the characteristic
downslope velocity ($m.s^{-1}$). The parameter $\hat{\nu}$ is the characteristic kinematic
viscosity:
\begin{equation}
  \label{eq:defnu}
  \hat{\nu} = \frac{\hat{\kappa}}{\hat{\rho}} \left ( \frac{ \hat{U}}{\hat{H}_0} \right )^{n-1},
\end{equation}
where $\hat{\kappa}$ ($kg.m^{-1}.s^{n-2}$ or $Pa.s^n$) is the consistency associated to
the Herschel-Bulkley rheology.

The typical height of the cavity $\hat{H}_0$ can be computed through:
\begin{equation}
  \label{eq:defH0}
  \hat{H}_0 = {\left ( \frac{\hat{\kappa} \hat{L}}{\hat{\rho} \hat{g} \cos \phi} {\left ( \frac{\hat{Q_0}}{\hat{W}} \right )}^n \right )}^{\frac 1 {2(1+n)}},
\end{equation}
where $\hat{g} = 9.81$ ($m^{2}.s^{-1}$) is the gravitational acceleration (note that
$\hat{g} = \|\overrightarrow{g}\|$ of figure \ref{fig:geom}), $\hat{W}$ ($m$) is the
transverse width of the cavity and $\hat{Q_0}$ ($m^3.s^{-1}$) is the known injected flow
rate of material. Note that if one starts from the given dimensional input data, one can
begin with \eqref{eq:defH0} and using:
\begin{equation}
  \label{eq:defUfQ}
  \hat{U} = \frac{\hat{Q_0}}{\hat{W} \hat{H_0}}
\end{equation}
one can compute subsequently the non dimensional quantities of the problem.

\bigskip

To fix the ideas, we recall the range of values for the parameters $(B,S,n)$. In the
context of CPB a "secured choice" (meaning extending a bit the range of potential extreme
values) of values is $B \in [0.5, 250]$ and $S \in [\sim 0.1 \mbox{ or } 1,
120]$. Concerning the power law index, it is generally observed $n \approx 1$. More
commonly, one might find CPB parameters of real mines close to $(B,S,n) \sim (150, 12, 1)$ (see
\cite{VSF21} for more details). We will come back to this point at the end of the article.

\subsection{Numerical solver}\label{sec:numsolv}

The solution of \eqref{eq:lub1D}-\eqref{eq:qdef} is approximated by a numerical solver.
It should be noted that this model has been extensively tested in various configurations
by leading groups in the field \cite{BCRS06,LBHH16,HoMa09}. So, even if it seems difficult
to find a rigorous proof of the existence and uniqueness of this model in the literature,
it does seem to be well-posed (in Hadamard's sense) in the configuration studied here.
The time evolution is treated by an explicit scheme (forward Euler). The spatial
discretization is treated using a centered scheme for $\partial h/\partial x$ and an
upwind scheme for $\partial q/\partial x$. Using the boundary conditions, $q$ can be
computed for all the points of the mesh, which enables to compute $h$ for all points
except the first one at $x=0$ (which would require the evaluation of $q(-\Delta x)$). This
point is treated by using the boundary condition $q(0)=1$ and a downwind scheme for
$\partial h/\partial x$. The resulting equation on $h(t,0)$ cannot be computed in closed
form (see \eqref{eq:qdef}) and leads to a scalar non linear equation whose root can be
solved with standard built-in libraries (e.g. \texttt{fzero} in Matlab or \texttt{brenth}
in Python-Scipy). \emph{Cf.} appendix \ref{sec:comph0} for more details.

For a given spatial discretization ($\Delta x$) we adapt the time step ($\Delta t$)
dynamically at each time iteration. A stability condition has been heuristically derived
so that the time step is determined by the spatial step. This is done by taking the
minimum of the two classical numerical stability constraints on $\Delta t$ associated to
the non-linear advection-diffusion problem \eqref{eq:lub1D}-\eqref{eq:qdef}. Namely, a
constraint due to the advection component of \eqref{eq:lub1D} (CFL condition) and a
constraint due to the explicit treatment of the diffusion component ($\Delta t$
constrained by $\Delta x^2$) \citep{Leveque2004}. \emph{Cf.} appendix \ref{sec:nsc} for
more details. It should also be noted that there is no difficulty in dealing with the
front where $h$ vanishes, due to the nature of this PDE which includes the threshold
\eqref{eq:ypp}.

\begin{figure}[!h]
  \begin{center}
    \includegraphics[scale=0.49]{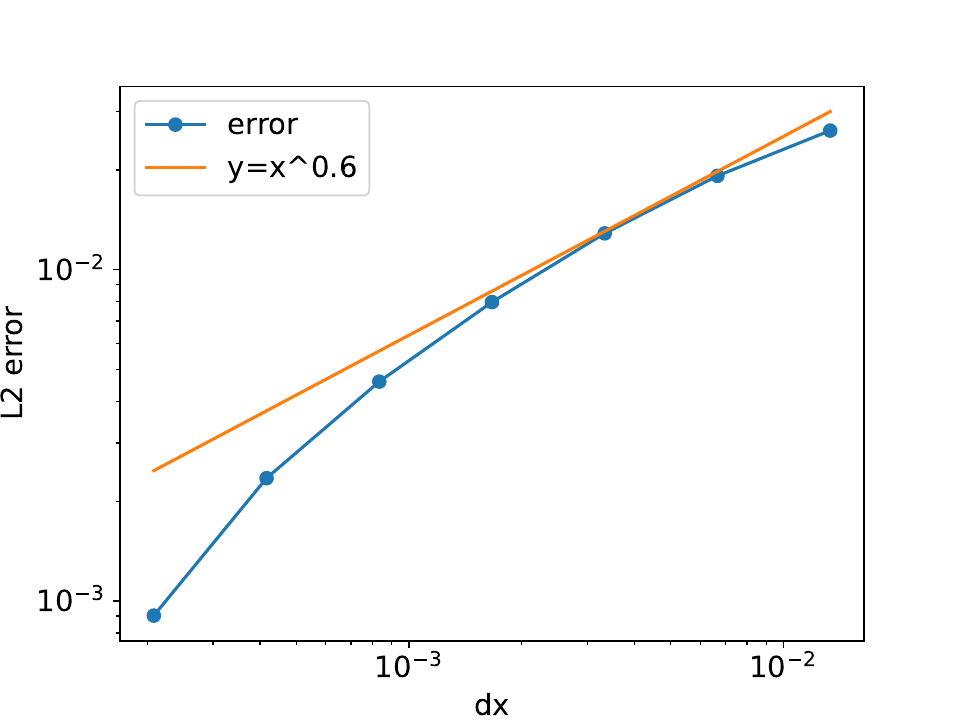}
  \end{center}
  \caption{Study of the convergence ($\Delta x$ is refined, and therefore also
      $\Delta t$ by \eqref{eq:stabcond}) of the PDE numerical solver for $B=100$,
    $S=120$, $n=0.8$. All solutions have been computed at their respective wall-touch
    time.  The reference solution ($h^{ref}$) has been computed using $n_x=9601$
    points. The $L^2$ error ($\|h^{ref}(x)-h^{\Delta x,\Delta t}(x)\|_2$) decreases at least like
    $y=x^{0.6}$.}%
  \label{fig:cv_100_120_0.8_all}
\end{figure}

To determine the spatial step, a numerical study of the convergence ($\Delta x$ is
refined, and therefore also $\Delta t$ by the stability condition \eqref{eq:stabcond}) was
performed. A reference solution has been computed beforehand using a very small $\Delta x$
(with $n_x=9601$ points in the mesh). Then, this solution is compared with the use of
others bigger $\Delta x$ ($n_x=76, 151, 301, \ldots, 4801$), to evaluate if there is a
numerical convergence of the norm of the difference of the two solutions, as $\Delta x$ is
reduced. The comparison is done at a final meaningful time: the one at which the fluid has
reached the right wall, a time referred to as \textit{wall-touch}. An example of such
study is illustrated on the figure \ref{fig:cv_100_120_0.8_all} for $B=100$, $S=120$ and
$n=0.8$. We see that the numerical scheme is convergent with an exponent power of at least
$0.6$. This exponent is not very high, which is known to be associated to the stiffness of
the non-linear problem. The development of higher order schemes is out of the scope of
this paper. Note that we also perform the same kind of study for ($B=30$, $S=15$, $n=0.8$)
and ($B=70$, $S=0.3$, $n=0.8$) which also show that the scheme is convergent with an order
between $0.6$ and $1$. This validates the use of this whole solver.

Some $(B,S)$ parameters lead to very fast resolutions, but some others can take up to two
hours for $n_x=301$ or multiple weeks for higher values of $n_x$. These correspond to high
values of $B$ and low values of $S$, i.e., a high yield stress and a low slope. For
such parameters, the wall-touch happens at higher time so many more time steps have to be
performed. As a compromise deduced from the above study, the use of $n_x=301$ 
provides very good accuracy while leading to reasonable computation times, in the
perspective of building a metamodel, as explained below.

This said, such numerical resolution remains costly, especially for some $(B,S,n)$
parameters. Assume now that we want to use this solver to estimate the parameters
$(B,S,n)$, given measured data of $h(x,t)$ profiles. We will then need to use inverse
problem algorithms to perform parameters' estimation. This may involve a lot of
evaluations of the direct problem: knowing $(B,S,n)$, compute $h(x,t)$. The computation
time may be too expensive. This is the justification to build a metamodel (or surrogate
model) which is an approximation of the PDE solver which is must faster to compute.

\subsection{Final direct problem}\label{sec:fdp}

With the previous considerations, the forward mapping provides, given $B$, $S$ and $n$,
the result of the numerical solver at wall-touch, using $n_x=301$ points for
the space discretization. This means that the output is a vector of $\mathbb{R}^{301}$ which
corresponds to an approximated height profile. The metamodel and inversions are performed
for fixed $n$, as it is assumed that physicists can provide a precise value for $n$, but
the methodology can be applied for any value of $n$. Note also that the model
\eqref{eq:lub1D} varies smoothly with $n$ for the values involved in the applications
($n \sim 1$). Thus the ultimate goal is to estimate $B$ and $S$ given height profile(s) of
a fluid in the cavity.

\section{Metamodel}\label{sec:metam}

A standard procedure for constructing a metamodel is the use of \emph{Polynomial Chaos Expansion}
(PCE), which consists in a polynomial approximation. They are often used for uncertainty
quantification \cite{lemaitre} or sensitivity analysis \cite{sudret_PCE}.
In order to reduce the dimension of the problem (which is currently $n_x=301$), we combine
this method with a \emph{Principal Component Analysis} (PCA).

\subsection{Polynomial chaos expansion (PCE)}

We start by describing the PCE procedure for some
bounded random vector $X$ of $\mathbb{R}^d$ (in our case
$X=(B,S)\in \mathbb{R}^2$) and a computational model
$\mathcal{M}:\mathbb{R}^d\rightarrow\mathbb{R}^m$.

We are interested in the random vector $Y=\mathcal{M}(X)$ (in our case, $m=301$ and 
$Y\in\mathbb{R}^{301}$ is the height profile at wall-touch).
We begin this description for the case $m=1$ in order to simplify the presentation
(section \ref{sec:geco} and first § of section \ref{sec:imco}); then we describe the full
case with $m\geq 1$ (last § of section \ref{sec:imco}).
Our goal is to approximate $Y$ by a formula of the type:
\begin{equation}
  Y = \sum_{\alpha\in \Lambda}\psi_{\alpha}(X)
\end{equation} 
where $(\psi_\alpha)_{\alpha\in \Lambda}$ are polynomials.  It is assumed that $Y$ has a
finite variance, which is reasonable for us since $(B,S)$ lives in a bounded domain and we
expect the mapping: $(B,S) \mapsto h$ to be smooth. We also assume that the coordinates of
$X$ are independent.

\subsubsection{General construction}\label{sec:geco}

In this section, we assume that $m=1$.
Let us define $f_{X_i}$ the marginal distribution of the $i^{th}$ component of $X$. We use
it to define an inner product for functions $u,v$ defined on the support $S_i$ of $X_i$
(we recall $X$ belongs to a bounded domain):
\begin{equation}
  \langle u,v \rangle_i=\int_{S_i}u(z)v(z)f_{X_i}(z)dz.
\end{equation}
Using this inner product as a scalar product, we define the orthogonality as usual. Then,
classic algebra procedures as Gram-Schmidt allows us to construct a family of orthogonal
polynomials $(P_k^i)_{k\in\mathbb{N}}$ which depends on $f_{X_i}$ 
and is not trivial in general. Still in some particular
cases, a known polynomial sequence is recovered. In our case (uniform distribution), we
recover the Legendre polynomials. For more examples of known families, see
\cite{sudret_PCE}, \cite{karniadakis}.
Note that originally, the methodology has been introduced with Hermite polynomials. For
other functions, we should talk about generalized PCE.

For now, we have constructed families of polynomials in one variable (one family for each
component of $X$). To construct polynomials of the whole $X$ vector, we will use products
of univariate polynomials (which makes sense as the components of $X$ are independent).
Formally, we introduce multi-indices $\alpha\in\mathbb{N}^{d}$ as
$\alpha=(\alpha_1,..,\alpha_d)$ and define:
\begin{equation} 
  \forall x\in\mathbb{R}^d,\ \psi_{\alpha}(x)=\prod_{i=1}^{d}P_{\alpha_i}^{i}(x_i).
\end{equation}
In this framework, it can be proven that we can write (see, e.g., \cite{lemaitre}):
\begin{equation}
  Y=\sum_{\alpha\in\mathbb{N}^{d}}c_{\alpha}\psi_{\alpha}(X),
  \label{PCE_theory}
\end{equation}
where the $c_{\alpha}$ are real and need to be determined.

\subsubsection{Implementation considerations}\label{sec:imco}

\paragraph{Standardized entries.} In fact, the Legendre polynomials (as well as other known
families) correspond to a precise distribution (\cite{karniadakis}). When the input $X$
does not correspond to one of these known distributions, some methods exist to determine
the corresponding family (see \cite{sudret_PCE} for references), but the standard
procedure consists in transforming $X$ into a new variable $\tilde{X}$ that enters in a
known case. In the case of uniform distributions, the standard distribution is the uniform
distribution over $[-1,1]$. However, $B\in[0.5,250]$ and $S\in[0.05,120]$. Thus we must
change the $(B,S)$ coordinates into $(\tilde{B},\tilde{S})$ using the formula :
\begin{equation} 
  \tilde{B} = \frac{B-125.25}{125.25},
\end{equation}
\begin{equation} 
  \tilde{S} = \frac{S-60.025}{60.025}.
\end{equation}

\paragraph{Truncation order.} In practice, we cannot determine an infinite number of
coefficients and one has to choose which polynomials to take. For this, let us define the
degree of a multivariate polynomial $\psi_{\alpha}$ by $\sum_{i=1}^{d}\alpha_i$ (the sum
of the degrees of the univariate polynomials defining $\psi_{\alpha}$). Usually,
polynomials with high degrees are associated with high order of interactions between the
inputs, which are generally limited \cite{lemaitre}, \cite{sudret_PCE}. That is the reason
why generally, the truncation is done by selecting polynomials of degree lower than a
threshold $\beta$, usually between 3 and 5.

Thus, once an order of truncation $\beta$ has been determined, we can rewrite the problem of
\eqref{PCE_theory} as:
\begin{equation} 
  Y\simeq\sum_{q=1}^{l}c_q\psi_q(X). \label{PCE_pratique}
\end{equation}
where we have reordered the polynomials for clarity of notation. Note that $l$ is finite,
it depends on the order of truncation $\beta$, but $l$ and $\beta$ are different.

\paragraph{Determination of coefficients.} There exist multiple techniques to compute the
coefficients of \eqref{PCE_pratique}. A first distinction is made between
\textit{intrusive} and \textit{non-intrusive} methods. The intrusive methods such as
Galerkin projection are based on the resolution of modified problems (i.e., not the
evaluation of $\mathcal{M}$ itself) which require the design of specific solvers, see
\cite{lemaitre}. More recently, non-intrusive methods have risen, based only on
evaluations of $\mathcal{M}$ and statistical tools. A popular solution is to treat the
problem as a least-square minimization problem \cite{migliorati_PCE},
\cite{PCE_sampling_review}, \cite{sudret_PCE}. Given a set of samples
$(X^{(1)},...,X^{(r)})$ with their evaluations $(Y^{(1)}, ...,Y^{(r)})$, we reduce the
problem to:
\begin{equation}
  \min_{c\in\mathbb{R}^{l}}\sum_{j=1}^{r}\left(Y^{(j)}-\sum_{q=1}^{l}c_q\psi_q(X^{(j)})\right)^2.
\end{equation} 
More details can be found in \cite{migliorati_PCE}, \cite{PCE_sampling_review} and
\cite{sudret_PCE}.

\paragraph{Multivariate output.} (case $m\geq 1$). The last aspect remaining to treat is
the fact that in our case, the output is a vector and not a scalar. It seems that in most
of the published articles of the literature, the components of the output are all treated
separately. This means that we decompose
$\mathcal{M}(X)=(\mathcal{M}(X)_1,...,\mathcal{M}(X)_m)$ and perform a different PCE for
each $\mathcal{M}(X)_i$, resulting in $m$ different computations of parameters. See,
e.g., \cite{sudret_PCE}, \cite{garcia_PCE_multivar_SA} and \cite{sun_PCE_multivar_SA}.

\subsection{Principal component analysis (PCA)}

We hereby give a brief review of the PCA. More details can be found in \cite{pca}.  In
this section, we want to study the random vector $Y\in\mathbb{R}^m$, the mapping
$\mathcal{M}$ is completely forgotten. We start by defining $\varphi^0\in\mathbb{R}^m$ as:
\begin{equation} \label{fst_pc}
  \varphi^0=\underset{\varphi\in\mathbb{R}^m}{\arg\max}\ Var(\langle\varphi,Y\rangle).
\end{equation}
We define :
\begin{equation}
  \alpha^0=\langle\varphi^0,Y\rangle.
\end{equation}
Then $\alpha^0$ is called the first Principal Component (PC) of $Y$ and $\varphi^0$ is its
direction.
For $1\leq k\leq m-1$, we recursively define the $k^{th}$ PC in a similar manner:
\begin{equation}\label{other_pcs}
  \varphi^k=\underset{\alpha^0,\dots,\alpha^{k-1}, \langle\varphi,Y\rangle\text{ uncorrelated}}
  {\arg\max}\ Var(\langle\varphi,Y\rangle),
\end{equation}
and
\begin{equation}
  \alpha^k=\langle\varphi^k,Y\rangle.
\end{equation}
With these variables we introduce the matrix notation:
\begin{equation}
  \alpha=\varphi Y,
  \label{pca_fond}
\end{equation}
where $\alpha$ is the column vector of the $\alpha_k$ and $\varphi$ is the squared matrix
whose lines are the $\varphi^k$.
One can show that in fact $\varphi$ is invertible. It follows that we can write
$Y=\varphi^{-1}\alpha$, so that the study of $\alpha$ is equivalent to the study of $Y$.
However, $\alpha$ is defined so that most of its variations are contained in its first few
components. This allows to reduce the dimension of the problem by performing a truncation
on the first PCs:
\begin{equation}
  Y\approx\varphi^{-1}\begin{pmatrix}
    \alpha^0\\
    \vdots\\
    \alpha^p\\
    \bar\alpha^{p+1}\\
    \vdots\\
    \bar\alpha^{m-1}
  \end{pmatrix},
\end{equation}
where $\bar\alpha^k=\esp{}{\alpha^k}$ and $p \in \llbracket 1,m-1\rrbracket$.

In our case we do not have access to the theoretical distribution of $Y$. As done in
practice \cite{pca}, we will thus use a data set which is supposedly a good sampling of
$Y$ and use sample variances and expectations to compute the PCs.

\subsection{Final surrogate via PCE - PCA coupling}

We now go back to the PCE. Until now, the metamodel is given by:
\begin{equation}
 \mathcal{S}(X)=(\Phi^0(X),...,\Phi^{m-1}(X))\approx \mathcal{M}(X),
\end{equation}
where $\Phi^k$ is the result of the PCE performed on the mapping
$X\mapsto\mathcal{M}(X)_k=Y_k$.
We replace these $m$ different PCEs performed on the components of $Y$ by PCEs performed
on the first PCs. More precisely, we assume that the $p+1$ first PCs account for most of
the variance of $Y$.
Then for $0\leq k\leq p$, we denote by $\theta^k$ the result of the PCE performed on the
mapping $X\mapsto\alpha^k(X)$ (by virtue of \eqref{pca_fond}, if $Y$ is viewed as a
function of $X$ then $\alpha^k$ can also be viewed as such). We use it for our final metamodel:
\begin{equation}
 \hat{\mathcal{S}}(X)=\varphi^{-1}\begin{pmatrix}
    \theta^0(X)\\
    \vdots\\
    \theta^p(X)\\
    \bar\alpha^{p+1}\\
    \vdots\\
    \bar\alpha^{m-1}
  \end{pmatrix}\approx \mathcal{M}(X).
\end{equation}
We recall that $\varphi$ and the $\bar\alpha^k$ are deterministic quantities precomputed
during the PCA step.

Such a combination between PCA and PCE is not common in the literature. Up to the
authors' knowledge it has been first used by \cite{PCA_PCE_1} and more recently in
\cite{PCA_PCE_2}.

\section{Results}\label{sec:resul}

Of note, in this section, we work with $n=1$. However, the entire methodology and codes
can be used with other values of $n$, e.g., between $0.2$ and $1.2$. Moreover as mentioned
in paragraph \ref{sec:fdp}, the present lubrication model varies smoothly for the values
of $n$ considered here ($n \sim 1$). It also turns out that the experimental pastes
studied in the last subsection are such that $n=1$ \cite{VSF21}.

\subsection{Metamodel performances}

In order to compute the coefficients of the metamodel, a sampling grid in the $(B,S)$
domain has to be established. The exact forward mapping is then evaluated for each couple
of the grid. We will work with a small regular grid of $400$ couples ($20$ different
values of $B$ and $20$ values of $S$) and a huge regular grid of $6084$ couples ($78$ by
$78$). Note that the couples which are the longest to be evaluated seem to be the same for
all values of $n$, which means that we can use the existing computations to estimate the
duration of the computation of each solution, so that the evaluation of a new dataset can
be efficiently parallelized. In practice, the small regular grid can be computed in five
hours using seven processors.

Let us now define the precision of a metamodel. Apart from our two regular grids, a
validation set of approximately 6000 couples sampled uniformly at random on the $(B,S)$
domain has been precomputed. Given a metamodel $\bar{\mathcal{M}}$, we evaluate it on all
the couples of the validation set and compare them to the real outputs
$\mathcal{M}(B,S)$. For each couple $(B,S)$ of the validation set, we define the
\textit{reconstruction error} $\|\mathcal{M}(B,S)-\bar{\mathcal{M}}(B,S)\|_{2}$. Some
statistics over the validation set for different parameters of the metamodel are presented
in the table \ref{tab:metamodel_stats}.  We comment them below.

Multiple orders of truncation for the PCE have been tested, from $4$ to $30$. As expected,
the higher the order is, the more precise the metamodel. From order $4$ to $6$, the error
is divided by $3$. From $10$ to $12$ it is divided by $2$. At order $30$, for the biggest
grid, there is approximately $3$\% of error for most couples of the validation set. There
is however a limit to it since there is a limited number of samples, so the system becomes
undetermined for a degree high enough. This happens for the smallest grid, which cannot
be associated to an order higher than 15. Another practical limit comes from the
computation time of the coefficients. While for small orders like $4$ or $5$, the
coefficients are estimated in a few minutes, it can take multiple days to train for order
$30$.

To combine PCA and PCE, one has to determine a number of principal components (PC) to
keep. For instance if one uses only $3$ PCs (i.e., $p=2$ following the notation of the
previous section), then the quality of the metamodel reaches a plateau after order
$10$. Using $10$ PCs (i.e., $p=9$), the precision is almost not affected even at order
$30$, so that is what we will use. Note that using the PCA with $10$ PCs, we have to
perform $10$ one-dimensional PCEs. If we did not use the PCA, we would have to perform
$n_x=301$ PCEs. So the use of the PCA divides the computation time by $30$.

\begin{table}
  \centering
  \caption{Different statistics of the errors of reconstruction for different
    metamodels, where $\beta$ is the order of truncation for the PCE.}
  \begin{tabular}{ |c|c|c|c|c|c| } 
    \hline
    Data set & $\beta$& PCA & median &$3^{\text{rd}}$ quantile& max \\
    \hline\hline
    \multirow{6}{*}{78x78} &\multirow{2}{*}{4} &no&1&1.8&16.8\\
    \cline{3-6} & &yes&1&1.8&16.8\\
    \cline{2-6} &\multirow{2}{*}{15} &no&0.018&0.027&0.59\\
    \cline{3-6} & &yes&0.018&0.027&0.59\\
    \cline{2-6} &\multirow{2}{*}{30} &no&0.0011&0.002&0.082\\
    \cline{3-6} & &yes&0.0018&0.0024&0.082\\
    \hline\hline
    \multirow{4}{*}{20x20} &\multirow{2}{*}{4} &no&1.1&1.9&15.8\\
    \cline{3-6} & &yes&1.1&1.9&15.8\\
    \cline{2-6} &\multirow{2}{*}{15} &no&0.021&0.032&1.13\\
    \cline{3-6} & &yes&0.021&0.032&1.13\\
    \hline
  \end{tabular}
  \label{tab:metamodel_stats}
\end{table}

We point out the fact that no matter the parameters of the metamodel, there are always a
few couples which lead to a significant error of reconstruction.  Even the most precise
version has a maximal error of $10$\% which can be problematic for an inversion
procedure. However these extreme cases happen very infrequently and almost only close to
the boundary of the $(B,S)$ domain, for small values of $B$ and generally also small
values of $S$. Since the domain has been chosen so that for the applications, the values
should be quite far from the boundary, we may expect the error to be under the
$3^{\text{rd}}$ quantile.

Until now the focus has been on the precision of the prediction of the metamodel. However,
let us recall the reason why a metamodel is needed in the first place is that an inversion
procedure will typically require hundreds of evaluations of the model, so that it needs to
be evaluated very fast. The metamodels presented above are all a lot faster to evaluate
than the real model, but they are not equal. While being slightly more precise, the PCE of
order $30$ without PCA requires the evaluation of $301$ bivariate polynomials of degree
$30$.  In comparison, the PCE of order 15 with PCA only requires $10$ polynomials of order
$15$ and a matrix multiplication. In practice, the latter is three times faster than the
first one.

\begin{figure}[h!]
  \centering
  \includegraphics[width=0.46\textwidth]{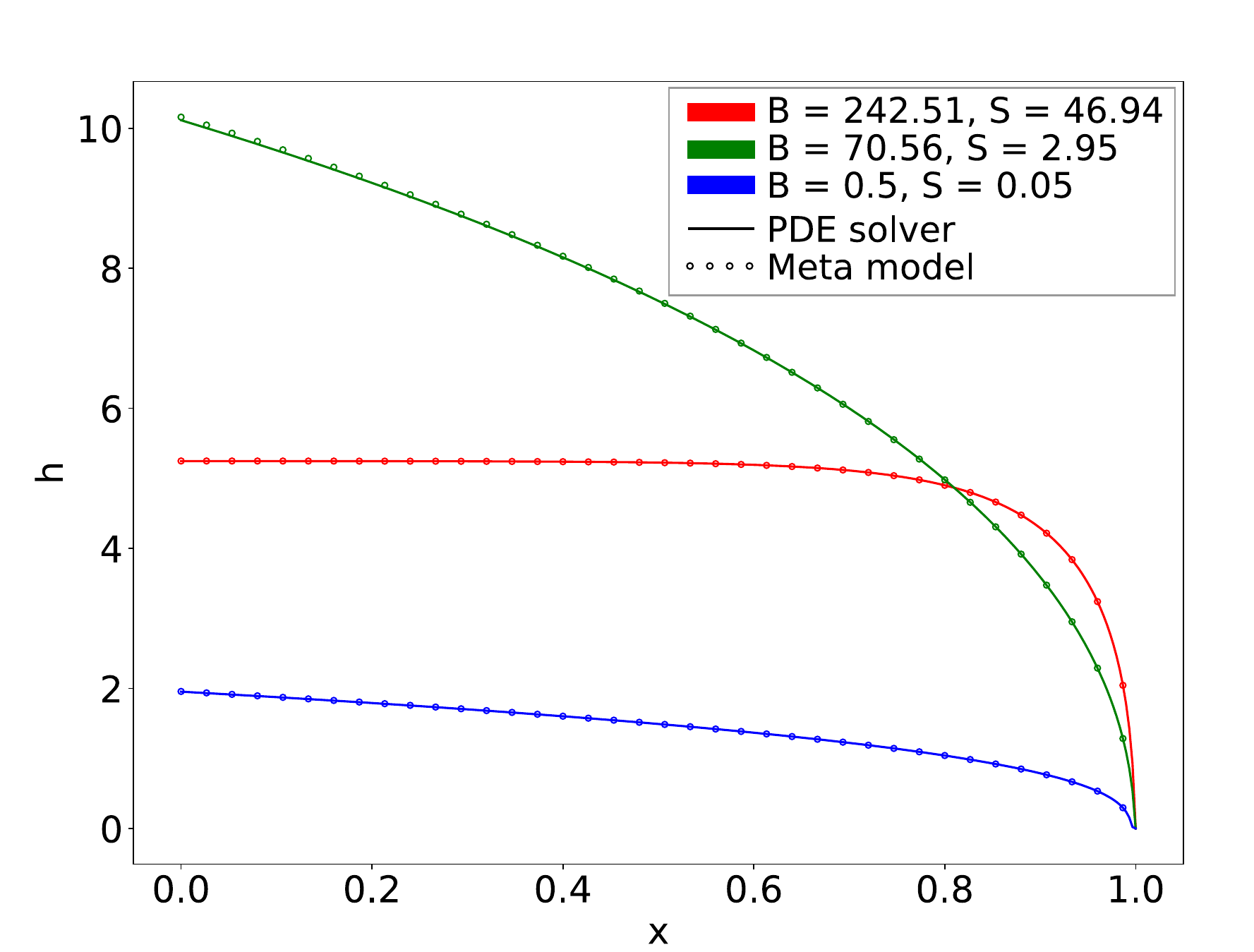}
  \caption{Comparison of the results from the PDE direct solver (solid line) and from the
    metamodel (circle markers), for three sets of $(B,S)$ parameters
    (colors).}\label{fig:compPDEmeta}
\end{figure}

In the figure \ref{fig:compPDEmeta}, we present typical reconstructions of the PDE
solution by the metamodel approach trained with the $20$ x $20$ data set, using the PCE of
order $15$ with PCA. Three representative $(B,S)$ couples are chosen and the associated
$h$ profiles given by the PDE solver and by the metamodel are superposed. We observe a
very good adequacy of both types of profiles.

Considering all of the above, we will now use the PCE-PCA with $\beta=15$ (following the
notation of the previous section) and $p=9$, trained on the smallest grid. This provides a
precise metamodel which can be evaluated very quickly and which can be trained in a day or
two at most (counting the dataset evaluation) on a domestic laptop.

\subsection{Parameters' estimation on synthetic data} \label{sec:synt}

The purpose of the metamodel is to be used for parameters' estimation. As a first step we
will test it on synthetic data, in the sense that we will use the profiles of the
validation set (termed as "real" in the following) and try to retrieve the values of $B$
and $S$, using our metamodel.

In addition to the real profiles $h(B,S)$, we create noisy profiles out of them. As field
data are expected to be available up to a precision of roughly $5$\% (see next section), we
will try $2$\%, $5$\%, and $10$\% of noise intensity. More precisely, we sample
$h_{noised}(B,S)=h(B,S)+\varepsilon$ where $\varepsilon\sim\mathcal{N}(0,\alpha h(B,S))$
and $\alpha=0.02$, $0.05$ or $0.1$. All in all, for each $(B,S)$ couple, we have at our
disposal $h(B,S)$ (the real one) and 3 noisy profiles with respectively $2$\%, $5$\%, and $10$\%
noise intensity. We will do so for $500$ couples of the original validation set.

The parameters' estimation is performed by the Nelder-Mead algorithm, with a maximum of
$400$ iterations, see \cite{nelder_mead} for details on the algorithm.  The error of the
estimation is defined by the Euclidean distance between the real $(B,S)$ couple and the
estimated one. The results are summed up in the table \ref{tab:invers_stats}.
The error displayed is the relative error of the estimated couples, i.e., 
\begin{equation}\label{eq:errBS}
err=\left\| \left(\frac{B_{real}-B_{estim}}{B_{real}},\frac{S_{real}-S_{estim}}{S_{real}}\right) \right\|_2.
\end{equation}

\begin{table}
  \centering
  \caption{Summary statistics of the errors (eq. \eqref{eq:errBS}) of parameters' estimation.}
  \begin{tabular}{ |c|c|c|c|c| } 
    \hline
    Noise &median &$3^{\text{rd}}$ quantile& max & var \\
    \hline
    0\% & 0.00045 & 0.00088 & 0.0299 & 7.46e-06\\
    2\% & 0.016 & 0.043 & 1.21 & 0.0276 \\
    5\% & 0.035 & 0.097 & 3.07 & 0.064  \\
    10\% & 0.067 & 0.18 & 5.03 & 0.202  \\
    \hline
  \end{tabular}
  \label{tab:invers_stats}
\end{table}

\begin{figure}[h!]
  \centering
  \includegraphics[scale=0.27]{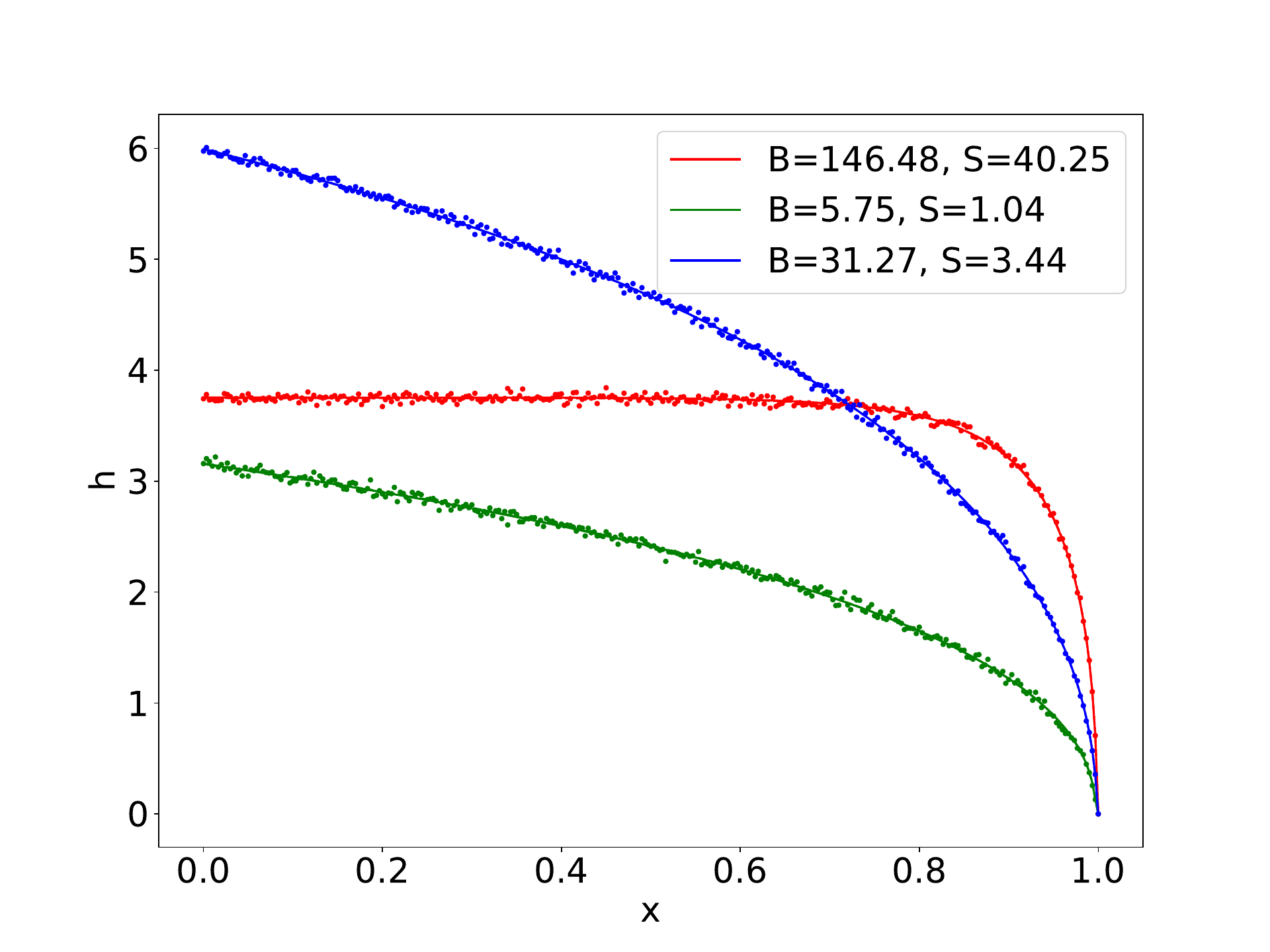}
  \caption{Comparison between a real height profile (solid line), its noisy version (circle) and
    the profile determined by the Nelder-Mead algorithm (dotted line). The noisy profiles
    are created using $2$\% of noise. The three colors correspond to three different $(B,S)$
    couples randomly chosen. Note that at this level of zoom the line and the dotted line
    are superimposed, so only one is visible.}
  \label{fig:inv_noise_2}
\end{figure}

\begin{figure}[h!]
  \centering
  \includegraphics[scale=0.27]{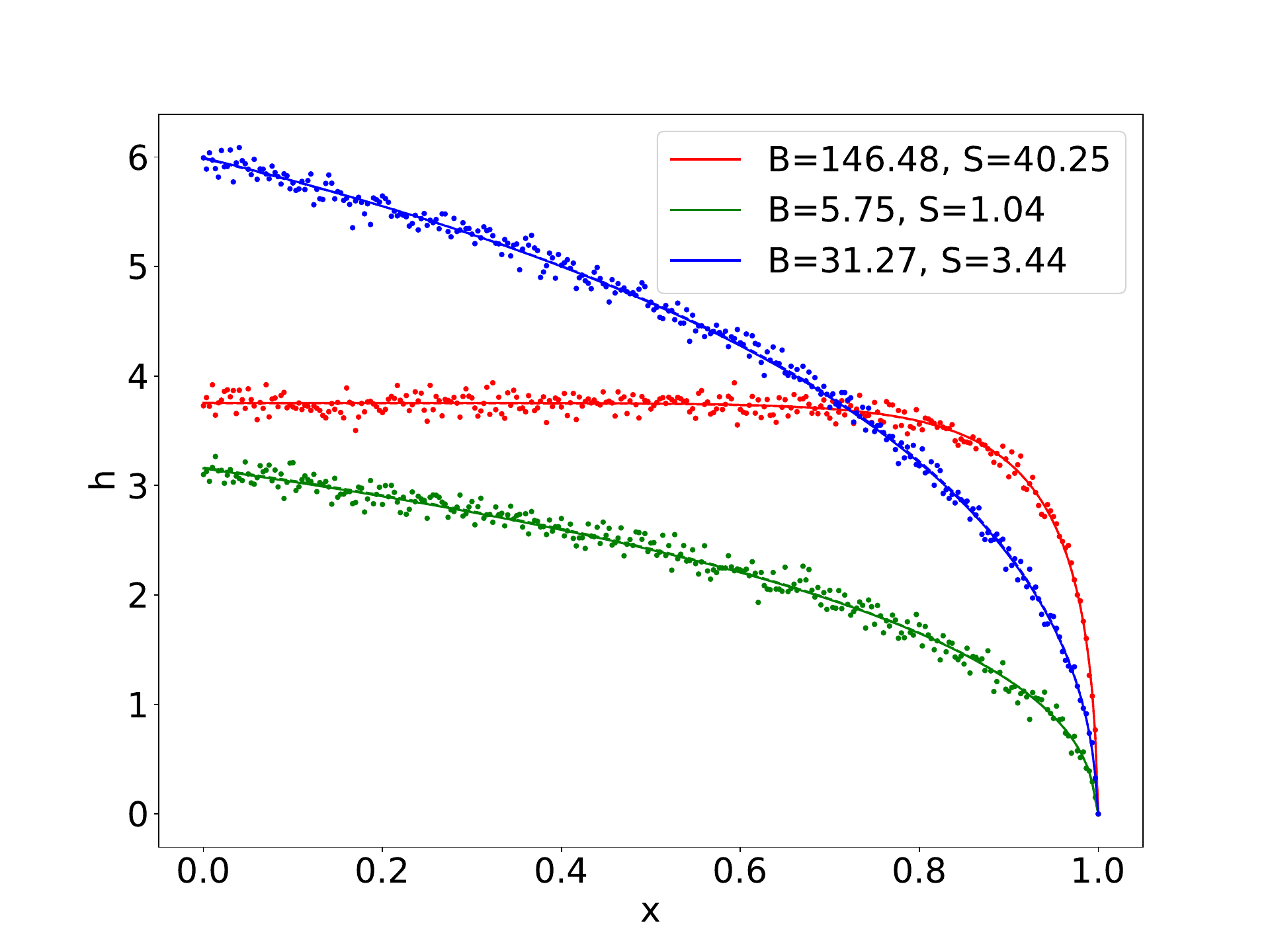}
  \caption{Comparison between a real height profile (solid line), its noisy version (circle) and
    the profile determined by the Nelder-Mead algorithm (dotted line). The noisy profiles
    are created using $5$\% of noise. The three colors correspond to three different $(B,S)$
    used in the figure \ref{fig:inv_noise_2}. Note that at this level of zoom the line and
    the dotted line are superimposed, so only one is visible.}
  \label{fig:inv_noise_5}
\end{figure}

For the non-noisy couples, the estimation is almost exact and for $2$\% of noise, there is
generally less than $5$\% of error. The situation for $5$\% of noise is a little worse, with
roughly $10$\% of error in the estimation. Note that these noisy profiles are extremely
chaotic compared with what we would expect from experimental data.  Still, the profiles
corresponding to the couples determined by the algorithm are very close to the real
profiles as displayed on the figures \ref{fig:inv_noise_2} and \ref{fig:inv_noise_5}.
If the $h$ profiles have $10$\% of noise, practitioners need to be cautious because the
estimation of $(B,S)$ can diverge more significantly than for the previous noise
intensities, even though real data would be expected to be less chaotic so the figure 
\ref{fig:inv_noise_10} (in the Appendix) should be thought of as a worst-case scenario.

\begin{figure}[h!]
  \begin{center}
    \centering
    \subfloat[][\centering paste 1
    ]{\label{fig:exp_1}\includegraphics[scale=0.4]{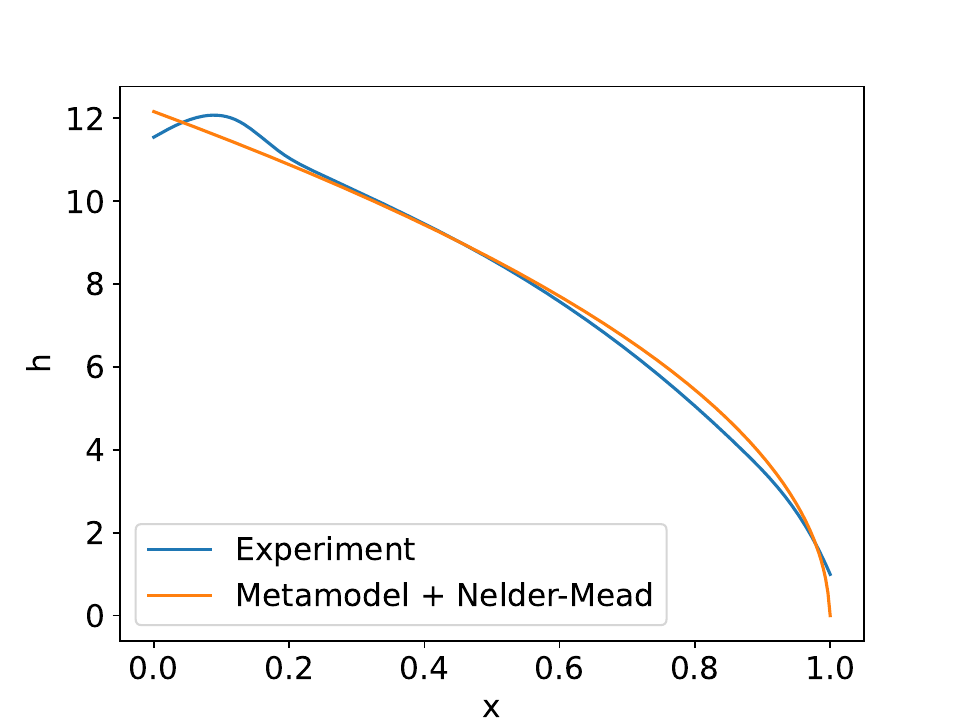}}%
    \qquad
    \subfloat[][\centering paste 2
    ]{\label{fig:exp_2}\includegraphics[scale=0.4]{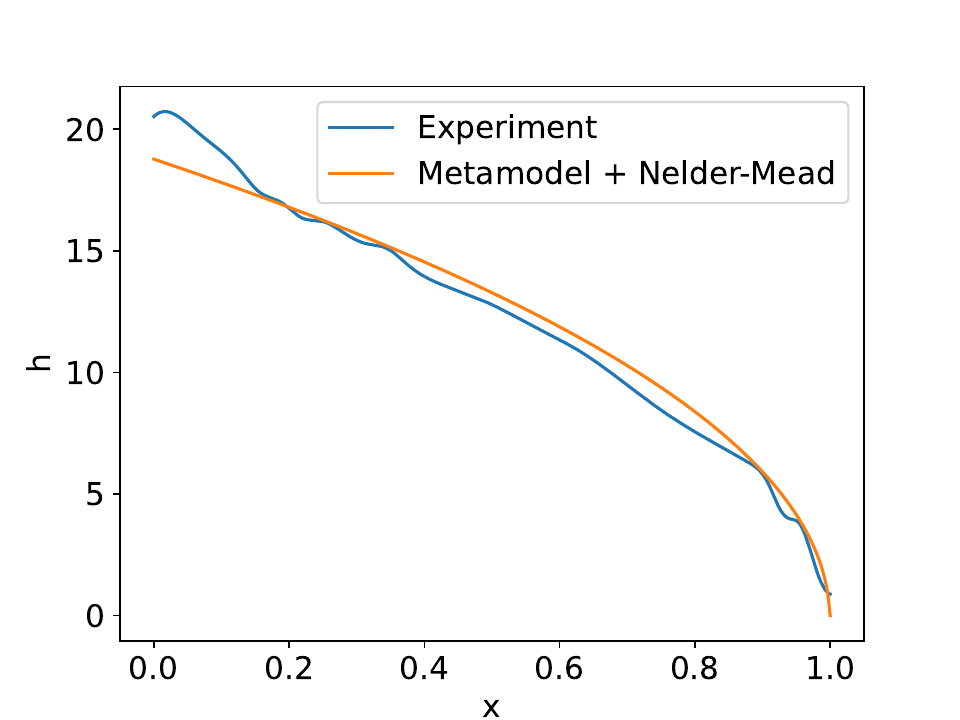}}%
    \qquad
    \subfloat[][\centering paste 3
    ]{\label{fig:exp_3}\includegraphics[scale=0.4]{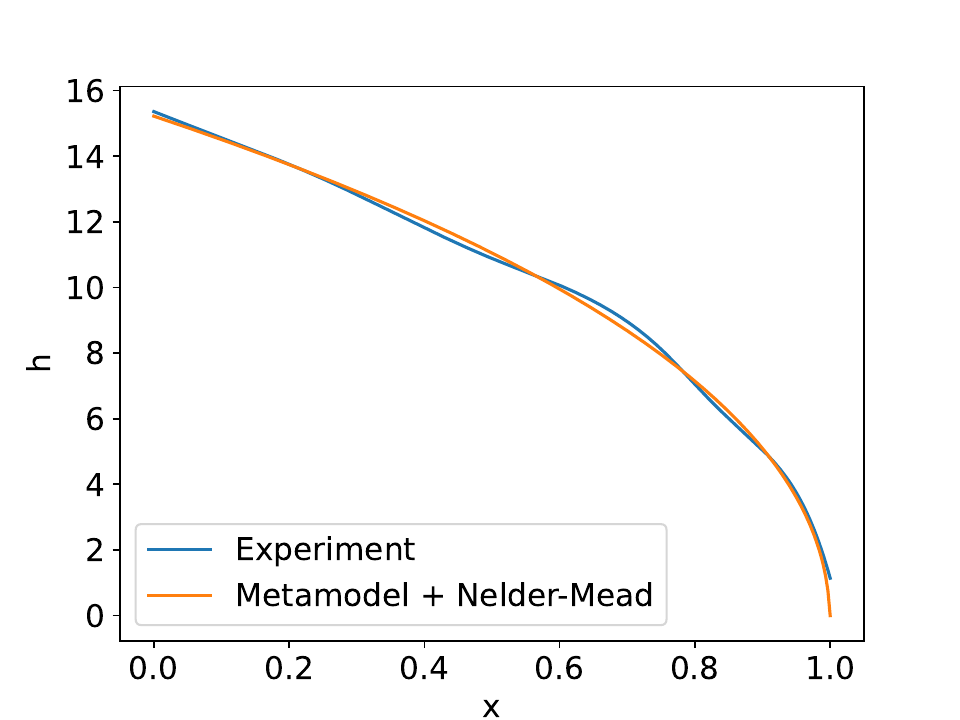}}%
    \qquad
    \subfloat[][\centering paste 4
    ]{\label{fig:exp_4}\includegraphics[scale=0.4]{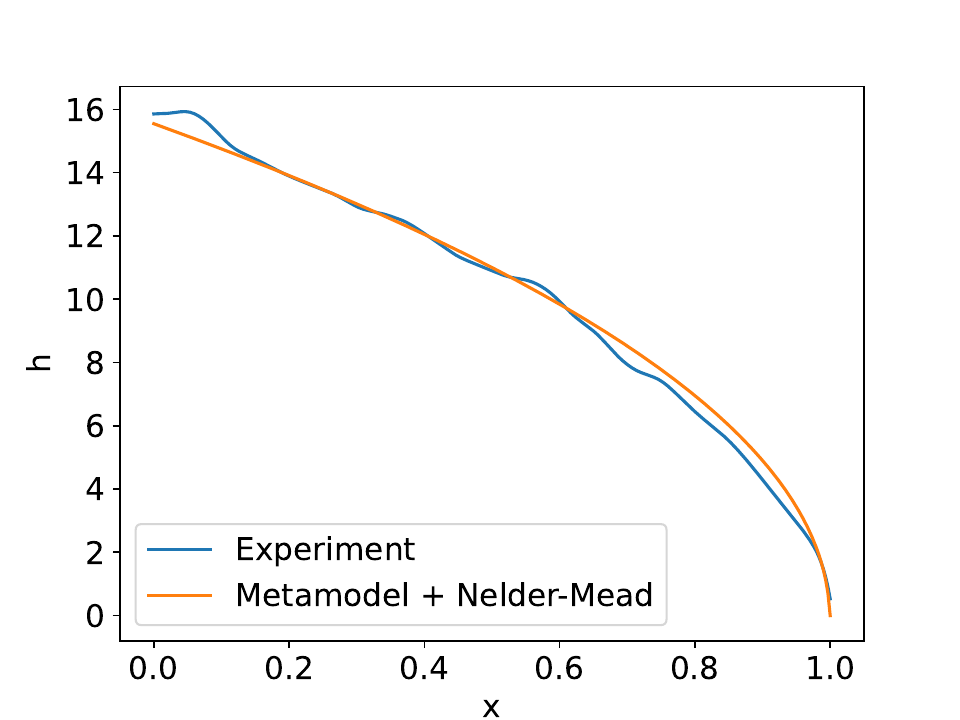}}%
    \caption{Fitting via the parameters' estimation procedure using four different real data 
      sets from \cite{VSF21}. For each paste are displayed in blue the experimental profile $h(x)$
      (non-dimensional variables, see section \ref{sec:descripmod}) and in
      orange the curve associated to the $(B,S)$ couple estimated by our procedure.}
    \label{fig:exp}
  \end{center}
\end{figure}

\subsection{Parameters' estimation on experimental data}

In this section, we apply the current algorithm to the experimental data measured in the
lab presented in \cite{VSF21}. In this article, four pastes are used to perform four
fillings of an empty box -- a reduced model in the lab for CPB experiment, with
pastes coming from the field. For each experiment the snapshot of the profile of $h$ at
wall-touch is obtained (see, e.g., figure 16 of \cite{VSF21}). This gives four curves
$x \mapsto h(x)$ that we can use to perform a $(B,S)$ parameters' estimation.

The results are presented in figure \ref{fig:exp}. As a first comment, we can note that
the experimental data do exhibit a sort of noise whose amplitude is of the same order as
in the virtual noise $2$\% - $5$\% used in section \ref{sec:synt}. The noise in the
experiments is not significantly associated to the image processing but is essentially due
to the very "thick" rheological behavior of the pastes (including mounding described
below). The parameters values estimated by the metamodel are given in table
\ref{tab:paramLab}.

\begin{table}
  \centering
  \caption{Parameters' estimation with $h$ profiles coming from lab data of
    \cite{VSF21}. "\#" stands for the Paste number. "exp" stands for the estimation done
    physically in the lab (with two techniques for $B$). "meta" stands for the estimation
    obtained with the present metamodel. "var" stands for the percentage of variation
    between \cite{VSF21} and present estimation with the metamodel. See main text for details.}
  \begin{tabular}{ |c|c|c|c|c|c|c| } 
    \hline
    \# & $B$: exp & meta & {\% var} & $S$: exp & meta & {\% var}\\
    \hline
    1 & [43;58]   & 73  & {69; 26} & 0.22 & 0.05 & {-77} \\
    2 & [102;122] & 176 & {73; 44} & 0.28 & 0.05 & {-82} \\
    3 & [74;98]   & 137 & {85; 40} & 0.28 & 2.2  & {686} \\
    4 & [68;93]   & 120 & {76; 29} & 0.26 & 0.05 & {-81} \\
    \hline
  \end{tabular}
  \label{tab:paramLab}
\end{table}


It is important to note that $(B,S)$ values were estimated in \cite{VSF21} with two
distinct techniques. Namely one estimation with the fitting of a Bingham constitutive law,
and the other with a creep shear experiment. Of note, the estimation of $B$ in the lab is
not trivial due to the use of real material as encountered in the field, for which very
precise rheometric measurements is a complicated task. The illustration of this can be
seen in the significant variability of $B$ shown for the two techniques in the second
column of Tab. \ref{tab:paramLab}. Another difficulty is the fact that the experiments may
present 3D mounding effect (see pastes 1, 2 and 4 in figure \ref{fig:exp}, around $x=0.1$)
which cannot be captured by the current 2D model (and not 3D) used here. In the present
study, this mounding effect translates in a "bump" shape perturbation for $x$ close to
$0.1$ which is equivalent to a noise on the curve of $h(x)$. This may challenge the
inversion procedure which tries to fit the PDE model (without mounding) to the
experimental curve.

It is also worth to mention that the lab experiment was such that $S$ is very small due to
a very small slope of the experimental device (in operational CPB context, one usually
encounters slightly bigger slopes and thus bigger $S$, as mentioned in section
\ref{sec:descripmod}). This implies that the numerical test for the inverse problem
is in a challenging part of the $(B,S)$ parameters' space because:

\begin{itemize}
\item (i) the smallest grid has been used for training, so the metamodel has been trained
  with very few values close to the boundary of the $(B,S)$ domain (the two smallest
  values of $S$ used for training are $0.05$ and roughly $6$, nothing in between), as a result
  it can be expected that the metamodel has only a limited precision in this very extreme
  range of S values;
\item (ii) a close inspection of the solutions $x\mapsto h(x)$ for different values of
  $(B,S)$ shows that for very small $S$ (particularly $S < 1$), only a slight change in
  the value of $S$ results in a significantly different solution of the PDE, so that both
  training the metamodel and solving the inverse problem become far more challenging in
  that area.
\end{itemize}
All the above reasons explain why this particular test is difficult. As well as the
significant percentage of variation seen in the Tab. \ref{tab:paramLab} between the
$(B,S)$ values of \cite{VSF21} and the values estimated by the present Metamodel+inversion
procedure.

This said, it should be noted that the curves obtained by the metamodel are very well
fitted with the experimental curves. It should also be noted that the inversion procedure
is very fast compared to an inversion using a resolution of the PDE. So, given the noise
that can be encountered with CPB field materials and the fact that their typical $(B,S)$
live in the bulk (and not at the boundaries) of the parameters' space of the
metamodel (giving more accuracy for exploring couples that minimize the discrepancy
between the curve of the model and the curve of the experiment), we believe that this
metamodel-inversion procedure can be one of the tools to make a first estimation of the
rheology of CPB experiments, based on the curve $x \mapsto h(x)$ at wall-touch.

For this purpose, we deliver on the Zenodo open repository \cite{zenodo2} the codes and
metamodel data \cite{BCVz} so that other teams having wall-touch snapshots of $h(x)$ can
test the current methodology with their own data.

\clearpage

\section{Summary and discussion}

Based on a partial differential equation (PDE) model derived from a lubrication
approximation, we have designed a metamodel approximating this PDE. The metamodel is
constructed using a combination of polynomial chaos expansion (PCE) and principal
component analysis (PCA). We have extensively tested different combinations of PCE and PCA
features, resulting in a metamodel that represents a good compromise between speed and
accuracy. As a result, given the two model parameters $(B,S)$, the solutions of the PDE
and the metamodel are essentially superimposed, while the computation of the metamodel is
much faster. Next, we tested this metamodel for parameters' estimation on synthetic data,
including noisy data. It is demonstrated, using statistically relevant samples, that with
the typical noise encountered in typical CPB applications, the inversion procedure leads
to good parameter estimations. The solutions of the inverse problem are such that the
estimate of $(B,S)$ is very good, as is the fit of the height profiles (data that can be
obtained from a laboratory snapshot). Finally, we also tested this inversion procedure on
a small sample of four available experimental data, in the worst-case scenario, i.e., with
a low-slope configuration that lies at the boundary of the parameters' space for CPB (and
thus also for the metamodel) with respect to $S$. The inversion procedure naturally
struggles to estimate these $(B,S)$ parameters. Given the high variability of $B$ values
obtained in \cite{VSF21}, an average $B$ error of $53\%$ can be considered quite
encouraging for such rough pastes and stiff numerical configuration.

In addition, given the good statistical results obtained on the synthetic data, it is
expected that with laboratory experiments conducted on CPB configurations mimicking the
field, i.e., with $S \sim 10$, inversion with the metamodel can be a valuable tool for
performing $(B,S)$ parameters' estimation. For testing purposes, we are providing the data
and metamodel codes \cite{BCVz} so that other teams can test them on their laboratory data. Not
only the ready-to-use metamodel, but also the code to recalculate this metamodel with
other parameters, are also provided in the repository, so that another training set of
$(B,S)$ and the corresponding metamodel can be computed offline in a few days with a
domestic laptop, providing a metamodel specialized in the $(B,S)$ range of a specific type
of CPB applications. The code is written in \texttt{Python} and is well-documented so that
a minimal investment is required to do so.

Future research directions would include a quantitative study of the influence of the
$\phi$ angle (involved in $S$), as well as a more systematic study of general slope
variations due to layers already deposited at the bottom. It is common practice in CPB to
carry out several injections on layers of cemented paste that have already been deposited
\cite{MHS2013}.


\appendix

\section{Numerical solver}

\subsection{Numerical computation of $h(t,0)$} \label{sec:comph0}

As stated in section \ref{sec:numsolv}, the numerical computation of $h(t,0)$ is not
straightforward. The time being fixed, let's write $h_1=h(t,\Delta x)$ (known) and
$h_0=h(t,0)$ (unknown). A downwind scheme for $\partial h/\partial x$ and the boundary
condition $q(0)=1$ together with the equation \eqref{eq:lub1D} lead to the following
problem on $h_0$:
\begin{equation}
\begin{split}
1=&\tilde\delta \frac{n}{(n+1)(2n+1)}Y_d(h_0)^{1+1/n}\times\\
&\left|S-\frac{h_1-h_0}{\Delta x}\right|^{1/n}
\left((2n+1)h_0-nY_d(h_0)\right)
\end{split}
\label{eq:h_0}
\end{equation}
with 
\begin{equation}
\tilde\delta=sgn\left(S-\frac{h_1-h_0}{\Delta x}\right),
\end{equation}
and
\begin{equation}
Y_d(h_0) = \max\left(h_0-\frac{B}{|S-\frac{h_1-h_0}{\Delta x}|},0\right).
\label{eq:Y_d}
\end{equation}

This problem is highly non-linear. However, it is constant equal to $0$ whenever
$Y_d(h_0)=0$. On the other part of the real line, it can be checked that the function of
$h_0$ defined by \eqref{eq:h_0} is actually strictly increasing and continuous, implying
the uniqueness of the solution of our problem. Because of the power law involved in
\eqref{eq:h_0}, an external solver is required to numerically solve the problem. In order
to accelerate its resolution, we provide to the solver the point after which $Y_d>0$ (so that it
does not get stuck in a constant region). Using \eqref{eq:Y_d} and $h_0\geq0$, we derive
the formula:
\begin{equation}
h_0>\frac{h_1-S\Delta x}{2}+\sqrt{\left(\frac{h_1-S\Delta x}{2}\right)^2+B\Delta x}.
\end{equation}

Note that in the special case $n=1$, the problem \eqref{eq:h_0} is equivalent to finding
the root of a polynomial of degree 7. However, in practice, it is not faster than using the
methodology presented above.

\subsection{Numerical stability condition} \label{sec:nsc}

We develop the expression \eqref{eq:lub1D} of the model, so as to obtain explicitly an
advection-diffusion form. Denoting $V(x)$ the transport coefficient and $D(x)$ the
diffusion coefficient, \eqref{eq:lub1D} is equivalent to:
\begin{equation}
  \label{eq:advdif}
  \frac{\partial }{\partial t} h(x,t) + V(x)\frac{\partial }{\partial x}  h(x,t)
  - D(x)\frac{\partial^2}{\partial x^2}  h(x,t) = 0
\end{equation}
with $\delta=\mbox{sgn}(S-h_x)$, 
\begin{equation}
  \label{eq:advec}
  V(x) = \delta Y^{1/n}|S-h_x|^{1/n}h,
\end{equation}
and
\begin{equation}
  \label{eq:diffc}
  \begin{split}
    & D(x)= \frac{-1}{(1+n)(1+2n)}Y^{1/n}|S-h_x|^{1/n-1}\times \\
    & \left(2nh\frac{B}{|S-h_x|}+(1+n)h^2+ 2n^2\frac{ B^2}{(S-h_x)^2}\right).
  \end{split}
\end{equation}
At each time iteration, as mentioned in the main text, the time step is determined as:
\begin{equation}
  \label{eq:stabcond}
  \Delta t = \min\left(\frac{\Delta x}{2\max_x(V(x))}, C_d\frac{(\Delta x)^2}{\max_x(D(x))}\right)
\end{equation}
where (the classical) $C_d = 0.5$ ensures stability for most simulations. One can lower
$C_d$, e.g.  $C_d = 0.05$ to stabilize the simulation if needed. Note that for the current
paper, where the solutions are not computed after the time of wall-touch, $C_d = 0.5$
leads to stable simulations for all studied parameters $(B,S)$.

\section{Noisy synthetic profiles}

The supplementary figure, mentioned in the main text, concerning the study of noisy
profiles with 10\% noise.

\begin{figure}[h!]
  \centering
  \includegraphics[scale=0.27]{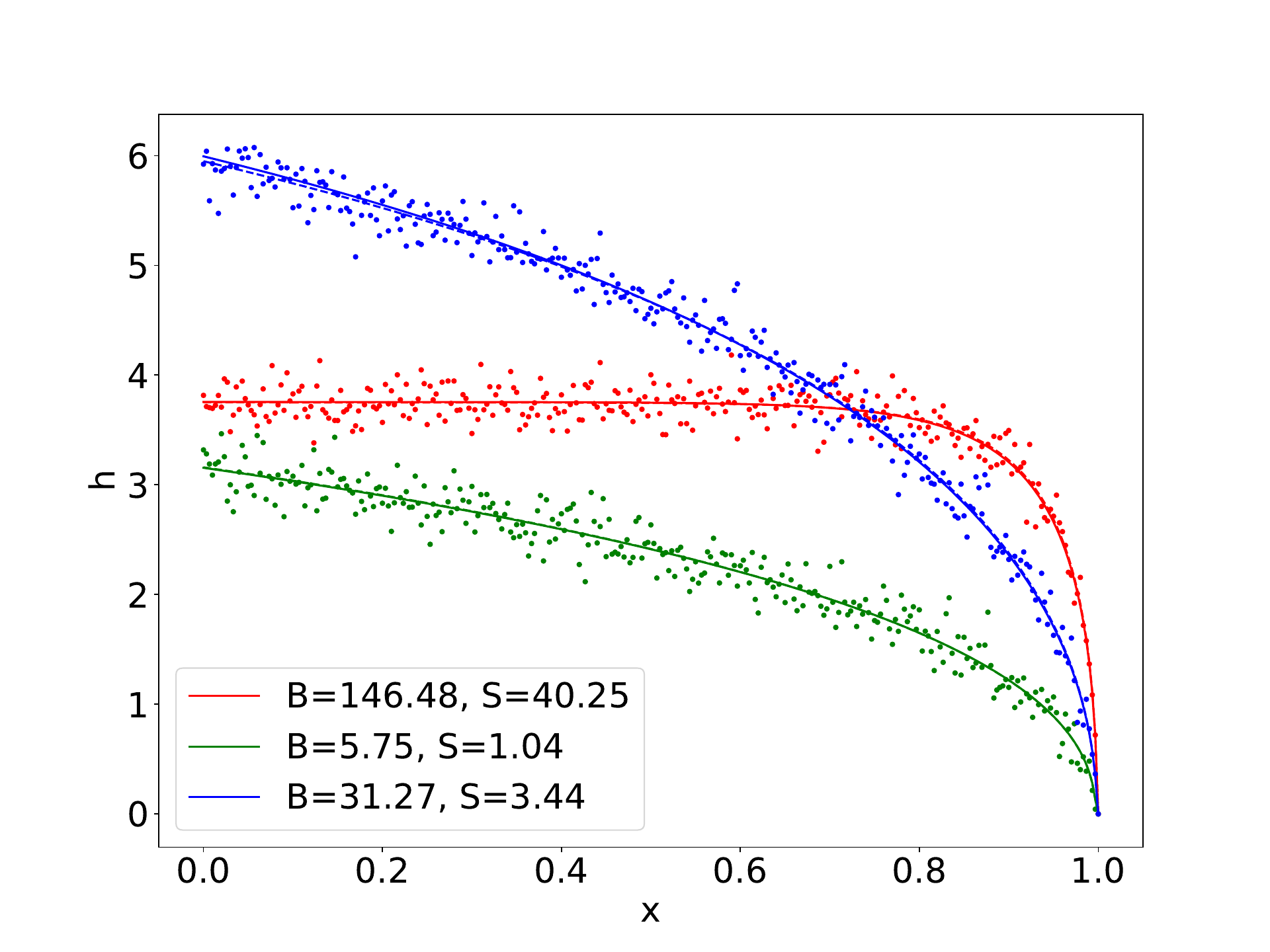}
  \caption{Comparison between a real height profile (solid line), its noisy version (circle) and
    the profile determined by the Nelder-Mead algorithm (dotted line). The noisy profiles
    are created using $10$\% of noise. The three colors correspond to three different $(B,S)$
    used in the figure \ref{fig:inv_noise_2}.}
  \label{fig:inv_noise_10}
\end{figure}

\section*{Acknowledgements}

We would like to extend our warmest thanks to Ian Frigaard and Yajian Shao for their
fruitful discussions on this work.

This study used components of the SciPy, NumPy, Scikit-learn, Chaospy, Numpoly and
Matplotlib Python libraries -- we thank the respective authors for making their code free
and open source.

We gratefully acknowledge support from the PSMN (Pôle Scientifique de Modélisation
Numérique) of the ENS de Lyon for the computing resources.
This research has been conducted with financial support from the French National Research
Agency (ANR) through the research project VPFlows, Grant number: ANR-20-CE46-0006.

\printcredits


\bibliographystyle{apalike}

\bibliography{cpbmet.bib}

\end{document}